\documentclass[11pt,reqno]{amsart}
\usepackage[T1]{fontenc}
\usepackage{lmodern,amsmath,amssymb,mathtools,booktabs}
\usepackage[margin=1.15in]{geometry}
\usepackage{listings}
\usepackage[hidelinks]{hyperref}
\allowdisplaybreaks[2]
\hypersetup{pdftitle={Finite-time blowup for six-dimensional wave maps into a negatively curved target},pdfauthor={Xuyi Wei},pdfkeywords={wave maps, negative curvature, finite-time blowup, exact polynomial certificates}}
\newtheorem{theorem}{Theorem}[section]
\newtheorem{proposition}[theorem]{Proposition}
\newtheorem{lemma}[theorem]{Lemma}
\newtheorem{corollary}[theorem]{Corollary}
\theoremstyle{remark}
\newtheorem{remark}[theorem]{Remark}
\newcommand{\R}{\mathbb R}
\newcommand{\Sph}{\mathbb S}
\newcommand{\dd}{\mathrm d}
\newcommand{\loc}{\mathrm{loc}}

\newcommand{\supp}{\operatorname{supp}}
\title[Blowup for six-dimensional wave maps]{Finite-time blowup for six-dimensional wave maps into a negatively curved target}
\author{Xuyi Wei}
\address{School of Mathematical Science, Beihang University, Beijing, China}
\email{weixuyi@buaa.edu.cn}
\date{}
\begin{document}
\begin{abstract}
We construct a complete rotationally symmetric metric on $\R^6$, with sectional curvature bounded between two strictly negative constants, for which corotational wave maps develop a finite-time singularity from smooth, compactly supported data. An explicit polynomial self-similar profile determines the target metric through an inverse construction. Negativity of curvature reduces to two polynomial inequalities on the unit interval, certified by Bernstein coefficients in exact rational arithmetic. Finite propagation and a continuation argument yield a full-space solution with maximal forward lifespan equal to one. A self-contained verification program is included.
\end{abstract}
\maketitle

\section{Introduction and main result}

Let $(N,h)$ be a Riemannian manifold. Equip $\R^{1+d}$ with the flat Minkowski metric of signature $(-,+,\ldots,+)$. The wave maps equation for a smooth map $U:\R^{1+d}\to N$ is
\begin{equation}\label{eq:wm}
 D_t\partial_t U-\sum_{j=1}^dD_j\partial_jU=0,
\end{equation}
where $D$ denotes the pullback of the Levi-Civita connection on $TN$. The conserved energy is
\begin{equation}\label{eq:energy-full}
 \mathcal E[U](t)=\frac12\int_{\R^d}
 \bigl(|\partial_tU|_h^2+|d_xU|_h^2\bigr)\,\dd x.
\end{equation}
The equation is invariant under the transformation $U(t,x)\mapsto U(t/\ell,x/\ell)$, $\ell>0$. Under this scaling the energy is multiplied by $\ell^{d-2}$. Thus the problem is energy critical in two spatial dimensions and energy supercritical in dimensions $d\geq3$; see \cite{TataruSurvey} for background.

For spherical targets, self-similar solutions provide an important mechanism of finite-time singularity formation in supercritical dimensions \cite{Shatah,BizonBiernat,Glogic2025}. The relationship between curvature and singularity formation is more subtle for negatively curved targets. Global regularity for a class of equivariant wave maps in dimension two was established in \cite{ShatahTZ}. In higher dimensions, the construction of Cazenave, Shatah and Tahvildar-Zadeh \cite{CST} showed that negative sectional curvature does not preclude blowup when $d=7$. Donninger and Glogi\'c \cite{DG} subsequently constructed negatively curved targets admitting explicit self-similar singularities for every $d\geq8$, together with a stability analysis. The complementary range $3\leq d\leq6$ was identified as open in the introduction to \cite{Glogic2025}.

The purpose of this paper is to treat the six-dimensional existence problem. Rather than prescribing the target metric and solving a nonlinear boundary-value problem for the similarity profile, we prescribe the profile and recover the metric. This inverse construction leaves a nontrivial positivity constraint: the reconstructed warping function must be strictly convex. We enforce that constraint by a finite polynomial certificate. In particular, neither a numerical approximation of a solution of the profile equation nor a spectral stability assumption is needed.

\begin{theorem}\label{thm:main}
There exist a complete smooth rotationally symmetric Riemannian metric $h$ on $\R^6$, a pole $p\in (\R^6,h)$, and constants $0<c\leq C<\infty$ such that
\begin{equation}\label{eq:pinch}
 -C\leq K_h(\Pi)\leq-c
\end{equation}
for every tangent two-plane $\Pi$. For this target there are smooth corotational initial data $(U_0,U_1)$, with $U_0=p$ and $U_1=0$ outside a compact set, whose wave map is smooth on $[0,1)\times\R^6$ and satisfies
\begin{equation}\label{eq:rate}
 |d_xU(t,0)|_h=\frac{\sqrt6}{1-t}.
\end{equation}
Consequently its maximal forward smooth lifespan is exactly one.
\end{theorem}

The target and the solution on the backward light cone are specified explicitly in Sections~\ref{sec:construction} and \ref{sec:completion}. The target is not assumed to have constant curvature. The theorem concerns existence, not stability under perturbation of the initial data.

\section{Corotational reduction and inverse construction}\label{sec:reduction}

\subsection{Rotationally symmetric geometry}

We work with warped product metrics of the form
\begin{equation}\label{eq:metric}
 h=\dd\rho^2+g(\rho)^2\,h_{\Sph^5},\qquad \rho>0,
\end{equation}
where $h_{\Sph^5}$ is the round unit sphere metric. We require $g$ to be a smooth odd function on $\R$, with
\begin{equation}\label{eq:pole}
 g(0)=0,\qquad g'(0)=1,\qquad g(\rho)>0\quad(\rho>0).
\end{equation}
These conditions ensure smoothness at the pole: writing $z=\rho\vartheta\in\R^6$ and $a(\rho)=g(\rho)/\rho$, one has, for $z\neq0$,
\[
 h_{ij}(z)=a(|z|)^2\delta_{ij}
   +\frac{1-a(|z|)^2}{|z|^2}z_i z_j.
\]
The function $a$ is smooth and even with $a(0)=1$. Both radial coefficients in this formula extend smoothly across $z=0$, and $h(0)$ is the Euclidean metric. We henceforth identify the pole with $z=0$.

For $\rho>0$ the sectional curvatures of radial and tangential two-planes, respectively, are (see \cite[Ch.~7]{ONeill})
\begin{equation}\label{eq:curvatures}
 K_{\mathrm{rad}}=-\frac{g''}{g},\qquad
 K_{\mathrm{tan}}=\frac{1-(g')^2}{g^2}.
\end{equation}
In the orthonormal splitting into radial and spherical directions the curvature operator is diagonal on the corresponding exterior-square subspaces. In particular the curvature of an arbitrary two-plane lies between these two values. It suffices to arrange $g''>0$ on $(0,\infty)$: since $g'(0)=1$, this implies $g'>1$. Curvature at the pole is determined by continuity.

\subsection{The scalar wave equation}

Writing $r=|x|$ and $\vartheta=x/r$, a corotational map is represented in the global Cartesian target coordinates by
\[
 U(t,x)=\psi(t,r)\vartheta\quad (x\neq0).
\]
Here the signed radial variable $\psi$ is allowed to be negative. The odd extension of $g$ makes the equations smooth under this convention. The map is smooth at $x=0$ whenever $\psi(t,r)=r\,b(t,r)$ with $b$ smooth and even in $r$. The angular identity map contributes five copies of the spherical metric, giving
\begin{equation}\label{eq:angular}
 |\partial_tU|_h^2=\psi_t^2,\qquad
 |d_xU|_h^2=\psi_r^2+\frac{5g(\psi)^2}{r^2}.
\end{equation}
The radial action, after omitting the area of $\Sph^5$, is
\[
 \frac12\iint \left(\psi_t^2-\psi_r^2-
       \frac{5g(\psi)^2}{r^2}\right)r^5\,\dd r\,\dd t.
\]
The Euler--Lagrange equation in this ansatz is
\begin{equation}\label{eq:radial}
 \psi_{tt}-\psi_{rr}-\frac5r\psi_r+\frac5{r^2}F(\psi)=0,
 \qquad F(\rho):=g(\rho)g'(\rho).
\end{equation}
This is also obtained directly from \eqref{eq:wm}: the angular part of the tension field vanishes for the identity angular map, and its radial part is exactly \eqref{eq:radial}. The signed coordinate version follows by smoothness. Thus a smooth solution of \eqref{eq:radial} satisfying the origin condition yields a wave map.

For smooth compactly supported radial data the conserved radial energy is
\begin{equation}\label{eq:radialenergy}
 E[\psi](t)=\frac12\int_0^\infty
 \left(\psi_t^2+\psi_r^2+\frac{5g(\psi)^2}{r^2}\right)r^5\,\dd r.
\end{equation}
Indeed, multiplying \eqref{eq:radial} by $r^5\psi_t$ and integrating gives $E'[\psi](t)=[r^5\psi_t\psi_r]_0^\infty=0$. The boundary terms vanish by smoothness at zero and finite propagation at infinity. The energy \eqref{eq:energy-full} is $|\Sph^5|E[\psi]$.

\subsection{Reconstructing the target from a similarity profile}

On the cone $0\leq t<T$, $0\leq r\leq T-t$, insert
\begin{equation}\label{eq:ansatz}
 \psi(t,r)=Q(y),\qquad y=\frac{r}{T-t}.
\end{equation}
The derivatives are
\[
 \psi_r=\frac{Q'}{T-t},\quad
 \psi_{rr}=\frac{Q''}{(T-t)^2},\quad
 \psi_t=\frac{yQ'}{T-t},\quad
 \psi_{tt}=\frac{y^2Q''+2yQ'}{(T-t)^2}.
\]
Thus \eqref{eq:radial} is equivalent to
\begin{equation}\label{eq:profile}
 (1-y^2)Q''+\left(\frac5y-2y\right)Q'-\frac5{y^2}F(Q)=0.
\end{equation}
Instead of solving this equation for $Q$ for a fixed $F$, choose an odd smooth profile with $Q(0)=0$, $Q'(0)=1$, and $Q'>0$ on $[0,1]$. Set
\begin{equation}\label{eq:inverse}
 F(Q(y))=\frac{1}{5}\bigl[y^2(1-y^2)Q''(y)+(5y-2y^3)Q'(y)\bigr].
\end{equation}
Since $Q$ is invertible there, this defines a single-valued smooth $F$ on its image. The required warping function must satisfy
\begin{equation}\label{eq:gintegral}
 g(\rho)^2=2\int_0^\rho F(s)\,\dd s.
\end{equation}
Making the substitution $s=Q(v)$ and integrating the term involving $Q''Q'$ by parts yields
\begin{align}
 g(Q(y))^2
 &=\frac25\int_0^y \bigl[v^2(1-v^2)Q''Q'
                         +(5v-2v^3)(Q')^2\bigr] \,\dd v\notag\\
 &=\frac15 y^2(1-y^2)Q'(y)^2
     +\frac85\int_0^y v Q'(v)^2\,\dd v.\label{eq:identity}
\end{align}
For the integration by parts we used $(Q'^2)'=2Q''Q'$ and $(v^2-v^4)'=2v-4v^3$. In particular the right-hand side is strictly positive for $0<y\leq1$. Smoothness and strict convexity of the reconstructed warping function remain to be verified.

\section{An explicit profile and its curvature certificate}\label{sec:construction}

\subsection{A polynomial criterion}

For a polynomial $P\in\mathbb Q[x]$ with $P(0)=1$, take
\begin{equation}\label{eq:Qdef}
 Q(y)=\int_0^yP(v^2)\,\dd v.
\end{equation}
Introduce the auxiliary polynomials
\begin{align}
 R(x)&=2x(1-x)P'(x)+(5-2x)P(x),\label{eq:R}\\
 S(x)&=(1-x)P(x)^2+\frac4x\int_0^xP(z)^2\,\dd z,\label{eq:S}\\
 C_0(x)&=(R(x)+2xR'(x))S(x)-P(x)R(x)^2.\label{eq:C}
\end{align}
The quotient in \eqref{eq:S} is polynomial, including at $x=0$. Since $R(0)=S(0)=5$, one has $C_0(0)=0$, and hence
\begin{equation}\label{eq:D}
 D(x)=C_0(x)/x\in\mathbb Q[x].
\end{equation}

\begin{lemma}[Polynomial curvature criterion]\label{lem:criterion}
Suppose that $P>0$ and $D>0$ on $[0,1]$. Then the reconstruction \eqref{eq:inverse}--\eqref{eq:gintegral} yields a smooth odd function $g$ on a neighborhood of $[-Q(1),Q(1)]$ such that $g'(0)=1$ and both sectional curvatures in \eqref{eq:curvatures} are strictly negative on $[0,Q(1)]$, including the pole. Moreover, the profile $Q$ solves \eqref{eq:profile}.
\end{lemma}
\begin{proof}
Since $Q'(y)=P(y^2)>0$, the odd map $Q$ has a smooth odd inverse on a neighborhood of the relevant interval. Substituting $x=y^2$ in \eqref{eq:inverse} and \eqref{eq:identity} gives
\begin{equation}\label{eq:parametric}
 F(Q(y))=\frac{yR(y^2)}5,\qquad
 g(Q(y))=y\sqrt{\frac{S(y^2)}5}.
\end{equation}
The positivity of $P$ and the formula for $S$ give $S(x)>0$ for $0\leq x\leq1$. Both inequalities persist on a slightly larger interval. Consequently the last formula, used with the odd inverse of $Q$, defines a smooth odd function $g$ also at the pole. It satisfies $g'(0)=1$ because $S(0)=5$ and $Q'(0)=1$. Differentiation of \eqref{eq:identity} verifies $gg'=F$.

For $y>0$, using \eqref{eq:parametric} and $F'=g'^2+gg''$ gives
\begin{align*}
 g'(Q(y))^2&=\frac{R(x)^2}{5S(x)},\\
 F'(Q(y))&=\frac{R(x)+2xR'(x)}{5P(x)},\\
 g(Q(y))g''(Q(y))&=\frac{xD(x)}{5P(x)S(x)}.
\end{align*}
In particular $g''>0$ away from the pole. More precisely,
\begin{equation}\label{eq:radialcurv}
 \frac{g''(Q(y))}{g(Q(y))}=\frac{D(x)}{P(x)S(x)^2}.
\end{equation}
The right-hand side extends positively and smoothly to $y=0$. Its limit is $D(0)/25$. The tangential curvature is negative since $g'>1$. Its limit at the pole agrees with the radial curvature: if $g(\rho)=\rho+a\rho^3+O(\rho^5)$ then both curvatures tend to $-6a$. This proves the geometric assertions. Equation \eqref{eq:profile} follows from \eqref{eq:inverse} and $F=gg'$.
\end{proof}

\subsection{The choice of polynomial}

We take the following polynomial of degree ten:
\begin{equation}\label{eq:Pexplicit}
\begin{split}
 P(x)={}&1+\frac7{10}x+\frac{63}{4}x^2
       -\frac{1356}{5}x^3+\frac{10836}{5}x^4
       -\frac{47502}{5}x^5\\
       &+\frac{248577}{10}x^6-\frac{199362}{5}x^7
       +\frac{192483}{5}x^8-\frac{102807}{5}x^9
       +\frac{18693}{4}x^{10}.
\end{split}
\end{equation}
Equivalently, a compact specification in the degree-ten Bernstein basis is
\begin{equation}\label{eq:Pbernstein}
 P(x)=\frac1{100}\sum_{j=0}^{10}c_j\binom{10}{j}x^j(1-x)^{10-j},
\end{equation}
where
\[
 (c_0,\ldots,c_{10})=(100,107,149,0,466,-385,844,-444,768,-95,680).
\]
The resulting similarity profile \eqref{eq:Qdef} is an odd polynomial of degree twenty-one. No optimization or floating-point search enters the following positivity proof; the polynomial is specified by the displayed rational coefficients.

\begin{proposition}[Exact positivity certificate]\label{prop:positivity}
For $P$ in \eqref{eq:Pexplicit} and the associated polynomials \eqref{eq:R}--\eqref{eq:D},
\begin{equation}\label{eq:pos}
 P(x)\geq1,\qquad D(x)>4\qquad(0\leq x\leq1).
\end{equation}
In addition, $D$ has degree $31$ and $D(0)=13$.
\end{proposition}
\begin{proof}
We use the elementary convex-hull property of the Bernstein basis (cf. \cite{Farouki}). For a polynomial $p$ of degree $n$, write its restriction to a subinterval $[a,b]$ as
\[
 p(a+(b-a)t)=\sum_{i=0}^n \beta_i\binom ni t^i(1-t)^{n-i},
 \qquad 0\leq t\leq1.
\]
The basis functions are nonnegative and sum to one, so
\begin{equation}\label{eq:convex}
 p(a+(b-a)t)\geq\min_i\beta_i.
\end{equation}
To compute these coefficients exactly, let $p(x)=\sum_{j=0}^n p_jx^j$. The monomial coefficients of $p(a+(b-a)t)$ are
\begin{equation}\label{eq:power}
 q_k=(b-a)^k\sum_{j=k}^n\binom jk p_ja^{j-k}.
\end{equation}
Using the identity
\[
 t^k=\sum_{i=k}^n\frac{\binom ik}{\binom nk}
                 \binom ni t^i(1-t)^{n-i},
\]
we obtain
\begin{equation}\label{eq:berformula}
 \beta_i=\sum_{k=0}^i\frac{\binom ik}{\binom nk}q_k.
\end{equation}
All entries in \eqref{eq:power}--\eqref{eq:berformula} are rational when $a,b,p_j$ are rational. Applying these formulas to $P$ and $D$ gives the certificate in Table~\ref{tab:cert}. The table records exact lower bounds rather than rounded numerical evaluations. In the $P$ rows the minima are exact; in the $D$ rows every Bernstein coefficient is strictly greater than the displayed integer.
\begin{table}[htbp]
\centering
\caption{Exact Bernstein coefficient certificates on subintervals of $[0,1]$.}\label{tab:cert}
\begin{tabular}{@{}ccc@{}}
\toprule
Polynomial & Subinterval & Bound for all Bernstein coefficients\\
\midrule
$P$ & $[0,1/2]$ & $\geq1$\\
$P$ & $[1/2,1]$ & $\geq33393/20480$\\
\midrule
$D$ & $[0,1/4]$ & $>4$\\
$D$ & $[1/4,1/2]$ & $>13$\\
$D$ & $[1/2,3/4]$ & $>48$\\
$D$ & $[3/4,1]$ & $>100$\\
\bottomrule
\end{tabular}
\end{table}
The arithmetic producing the table, the degree assertion and the value at the origin is given in Appendix~\ref{app:code}; it consists entirely of polynomial additions, multiplications, derivatives and exact rational comparisons. For instance, the first coefficient at zero follows from $P(x)=1+\frac7{10}x+O(x^2)$, giving $R(x)=5+\frac{29}{10}x+O(x^2)$ and $S(x)=5+\frac{16}{5}x+O(x^2)$, so that $C_0(x)=13x+O(x^2)$. The upper degree is at most $31$ after division by $x$, and the nonzero leading coefficient is checked in the same arithmetic. Finally \eqref{eq:convex} on the indicated subintervals yields \eqref{eq:pos}.
\end{proof}

\begin{remark}
The strictly positive margins in Table~\ref{tab:cert} are important: they guarantee curvature negativity at the lateral cone boundary itself, so that the reconstructed metric extends beyond the image of the closed similarity interval. Positivity at a discrete set of sample points would not suffice for this purpose.
\end{remark}

\section{Completion of the target and construction of blowup}\label{sec:completion}

\subsection{Extension to a complete negatively curved metric}

Let $L=Q(1)$ for the polynomial profile constructed above. Lemma~\ref{lem:criterion} and Proposition~\ref{prop:positivity} produce a smooth odd local warping function, denoted $g_{\loc}$, on a neighborhood of $[-L,L]$. Its Taylor expansion follows directly from \eqref{eq:parametric}:
\begin{align*}
 Q(y)&=y+\frac7{30}y^3+O(y^5),\\
 g_{\loc}(Q(y))&=y+\frac8{25}y^3+O(y^5),\\
 g_{\loc}(\rho)&=\rho+\frac{13}{150}\rho^3+O(\rho^5).
\end{align*}
In particular $K_h(0)=-13/25$.

Define the local even function
\[
 k_{\loc}(\rho)=\frac{g_{\loc}''(\rho)}{g_{\loc}(\rho)}
\]
with its removable value $k_{\loc}(0)=13/25$. By \eqref{eq:radialcurv}, it is strictly positive on $[-L,L]$. Smoothness and compactness give $\delta>0$ such that $k_{\loc}$ is defined and positive for $|\rho|\leq L+2\delta$. Choose an even smooth cutoff $\zeta$ with $0\leq\zeta\leq1$, $\zeta=1$ on $[-L-\delta,L+\delta]$ and $\supp\zeta\subset(-L-2\delta,L+2\delta)$. Extend $k_{\loc}$ by
\begin{equation}\label{eq:kextend}
 k(\rho)=\zeta(\rho)k_{\loc}(\rho)+1-\zeta(\rho).
\end{equation}
This is well-defined using the local function only on the support of $\zeta$. There are constants $0<c\leq C<\infty$ for which $c\leq k\leq C$ on $\R$. Solve
\begin{equation}\label{eq:gode}
 g''=kg,\qquad g(0)=0,\qquad g'(0)=1.
\end{equation}
Linear ODE theory gives a unique global smooth solution. It is odd by uniqueness and evenness of $k$, and it agrees with $g_{\loc}$ on $[-L-\delta,L+\delta]$. Initially $g>0$ on the positive axis. If it had a first positive zero, then on the interval before that zero we would have $g''=kg>0$ and hence $g'\geq1$, which is impossible. Thus $g>0$, $g''>0$, and $g'>1$ for all $\rho>0$.

The metric \eqref{eq:metric} therefore extends smoothly to all of $\R^6$. It is complete: every curve escaping compact subsets must have unbounded radial coordinate, and the length of any curve is bounded below by the total variation of its radial coordinate. Equivalently, the distance from the pole to $z$ equals $|z|$, so closed metric balls about the pole are compact, and Hopf--Rinow applies. The curvature pinching also follows directly. For $\rho>0$,
\begin{equation}\label{eq:pinchingproof}
 c g(\rho)^2\leq g'(\rho)^2-1
 =2\int_0^\rho k(s)g(s)g'(s)\,\dd s\leq C g(\rho)^2.
\end{equation}
Together with $g''/g=k$, the formulas \eqref{eq:curvatures} imply \eqref{eq:pinch} on every two-plane away from the pole. Continuity yields the same bounds at the pole. This completes the construction of the fixed target manifold.

\subsection{Smooth local solutions and a continuation criterion}

For clarity, we record the local theory needed to turn the cone solution into a full-space solution. It does not require a global regularity theorem for wave maps. Set $\psi=rv$ in \eqref{eq:radial}. A direct computation gives
\begin{equation}\label{eq:lift}
 v_{tt}-v_{rr}-\frac7r v_r+\frac5{r^3}\bigl(F(rv)-rv\bigr)=0.
\end{equation}
Since $F$ is smooth and odd with $F'(0)=1$, the function
\[
 J(s)=\begin{cases}(F(s)-s)/s^3,&s\neq0,\\ F'''(0)/6,&s=0\end{cases}
\]
is smooth and even. Thus \eqref{eq:lift} is precisely the radial restriction of the semilinear equation
\begin{equation}\label{eq:8wave}
 v_{tt}-\Delta_{\R^8}v+5v^3J(|X|v)=0,\qquad X\in\R^8.
\end{equation}
The nonlinearity is smooth even at $X=0$: a smooth even one-variable function is a smooth function of the square of its argument on the nonnegative axis. In particular $(X,v)\mapsto J(|X|v)$ is smooth. The nonlinearity vanishes at $v=0$.

\begin{lemma}[Continuation]\label{lem:continue}
Smooth, radial, compactly supported Cauchy data for \eqref{eq:8wave} generate a unique smooth radial solution on a maximal interval $[0,T_*)$. It has finite propagation speed. If $T_*<\infty$, then
\[
 \sup_{0\leq t<T_*}\|v(t)\|_{L^\infty(\R^8)}=\infty.
\]
\end{lemma}
\begin{proof}
We include the argument to account for the spatially dependent coefficient. On any fixed bounded time interval, finite propagation restricts the support to a fixed ball; one can multiply the nonlinearity by a smooth radial spatial cutoff equal to one on a slightly larger ball. Standard high-order semilinear wave local theory then applies (see \cite{Taylor,Sogge}). More explicitly, for an integer $m\geq6$, the wave energy estimate in $H^m\times H^{m-1}$, Sobolev embedding, and the local Lipschitz estimate for smooth composition give a contraction on a short time interval by Duhamel's formula. Local uniqueness and rotational invariance give radiality, while the local energy estimate on cones gives finite propagation. Higher regularity is propagated by differentiating the equation and using the tame composition estimates.

If $\|v(t)\|_\infty\leq B$ up to a finite time $T_*$ and the support is contained in a ball of radius $R_*$, the smoothness of the nonlinearity on $|X|\leq R_*+1$, $|v|\leq B$ gives
\[
 \|5v^3J(|X|v)\|_{H^{m-1}}\leq C_{m,B,R_*}\|v\|_{H^{m-1}}.
\]
This is the usual tame composition estimate for a smooth function of $(X,v)$ vanishing at $v=0$, after the spatial cutoff described above. Its constants depend on finitely many derivatives of that function on the compact set. Let
\[
 \mathcal H_m(t)=\|v(t)\|_{H^m}^2+\|v_t(t)\|_{H^{m-1}}^2.
\]
The linear inhomogeneous energy estimate (or differentiation of the equivalent energy including $\|v\|_{H^{m-1}}^2$) yields
\[
 \mathcal H_m(t)\leq C_m\mathcal H_m(0)\exp(C_{m,B,R_*}t).
\]
Thus the local existence time can be bounded from below in terms of a fixed $H^m\times H^{m-1}$ bound and the spatial cutoff. Restarting near $T_*$ extends the solution beyond $T_*$, a contradiction. Persistence of arbitrary higher regularity follows from the same estimate. This proves the lemma.
\end{proof}

The transformation $\psi=rv$ identifies a smooth radial solution of \eqref{eq:8wave} with a smooth corotational map in six dimensions: in Cartesian target coordinates $U(t,x)=v(t,|x|)x$ (here $x\in\R^6$). A smooth radial function in eight dimensions restricts to a smooth radial function on the six-dimensional coordinate subspace. The wave map equation then follows away from zero by \eqref{eq:lift} and at zero by continuity. Smooth uniqueness for the full wave maps system, if desired, follows in these global Cartesian coordinates by the usual high-order local theory for its semilinear derivative nonlinearity.

\subsection{Localized self-similar data and the singularity}

Set $T=1$ and define, on the backward cone,
\begin{equation}\label{eq:conesol}
 \psi_*(t,r)=Q\left(\frac{r}{1-t}\right),\qquad
 0\leq r\leq1-t,\quad 0\leq t<1.
\end{equation}
The range of $Q$ on the cone is contained in $[0,L]$, where $g$ agrees with the reconstructed warping function. Hence \eqref{eq:conesol} satisfies \eqref{eq:radial}. In fact the same solution is defined in a slightly larger similarity interval, by the local extension of the metric and the profile. There is no regularity loss at $r=1-t$. Since $Q$ is odd, $Q(y)/y$ is an even polynomial and the solution is also smooth at $r=0$ for $t<1$.

Choose a smooth even cutoff $\chi:\R\to[0,1]$ such that $\chi=1$ on $[-1,1]$ and $\chi=0$ for $|r|\geq2$. Prescribe the radial initial data
\begin{equation}\label{eq:data}
 \psi_0(r)=\chi(r)Q(r),\qquad \psi_1(r)=\chi(r)rQ'(r).
\end{equation}
Dividing by $r$ gives smooth even compactly supported data for \eqref{eq:8wave}, with the quotients at zero defined by continuity. The corresponding geometric data are $U_0(x)=\psi_0(r)x/r$ and $U_1(x)=\psi_1(r)x/r$ in Cartesian target coordinates. In particular $U_1(x)$ is a tangent vector at $U_0(x)$ under the coordinate identification of the tangent bundle. Both data are smooth, and $U_0=p$, $U_1=0$ outside $B(0,2)$.

Let $v$ be the maximal smooth solution supplied by Lemma~\ref{lem:continue} and write $\psi=rv$. By uniqueness and finite propagation, the cone formula holds throughout the intersection of the lifespan of $v$ and the cone:
\begin{equation}\label{eq:match}
 \psi(t,r)=\psi_*(t,r)\qquad(r+t\leq1,\ 0\leq t<\min\{1,T_*\}).
\end{equation}
One can apply finite propagation directly to the equation in eight dimensions, comparing $v$ with $(1-t)^{-1}Q(r/(1-t))/(r/(1-t))$, which is smooth at zero and up to the lateral cone boundary.

We now show that $T_*\geq1$, rather than merely concluding that a singularity occurs no later than one. Fix $t_0<1$ and set $r_0=(1-t_0)/2$. The following estimates hold as long as $t\leq t_0$ and $t<T_*$. For $r\leq r_0$, \eqref{eq:match} gives
\begin{equation}\label{eq:inner}
 |v(t,r)|\leq\frac1{1-t_0}\max_{0\leq y\leq1}\left|\frac{Q(y)}y\right|.
\end{equation}
For $r\geq r_0$, conservation of \eqref{eq:radialenergy}, compact support at each finite time, and Cauchy--Schwarz give
\begin{align*}
 |\psi(t,r)|
 &\leq\int_r^\infty|\psi_r(t,s)|\,\dd s\\
 &\leq \left(\int_r^\infty s^5|\psi_r(t,s)|^2\,\dd s\right)^{1/2}
         \left(\int_r^\infty s^{-5}\,\dd s\right)^{1/2}\\
 &\leq\frac{\sqrt{2E[\psi](0)}}{2r^2}.
\end{align*}
Consequently
\begin{equation}\label{eq:outer}
 |v(t,r)|\leq\frac{\sqrt{2E[\psi](0)}}{2r_0^3}\qquad(r\geq r_0).
\end{equation}
Equations \eqref{eq:inner} and \eqref{eq:outer} bound the full spatial supremum on the lifespan up to any fixed $t_0<1$. If $T_*<1$, choose $t_0\in(T_*,1)$, and Lemma~\ref{lem:continue} contradicts maximality. Thus the solution exists on $[0,1)$.

Finally $Q'(0)=1$. In the target Cartesian coordinates,
\[
 U(t,x)=\frac{x}{1-t}+O_t(|x|^3)\quad (x\to0),
\]
so $d_xU(t,0)=(1-t)^{-1}\mathrm{Id}_{\R^6}$. The target metric is Euclidean at the pole, yielding \eqref{eq:rate} with the Hilbert--Schmidt norm. This prevents any smooth continuation through $t=1$. The maximal forward lifespan is therefore exactly one, proving Theorem~\ref{thm:main}.

\begin{corollary}\label{cor:energy}
On the same fixed target there exist smooth compactly supported corotational initial data of arbitrarily small positive energy which develop a finite-time singularity.
\end{corollary}
\begin{proof}
If $U$ is the solution of Theorem~\ref{thm:main}, define $U_\ell(t,x)=U(t/\ell,x/\ell)$ for $\ell>0$. Its blowup time is $\ell$, and its initial energy is $\ell^4\mathcal E[U](0)$. Letting $\ell\downarrow0$ proves the assertion. This is consistent with small-data global regularity results formulated in the scaling-critical Sobolev topology, rather than in the supercritical energy topology.
\end{proof}

\section{Concluding remarks}

The inverse reconstruction \eqref{eq:inverse} separates the geometric curvature constraint from the wave equation. The certificate in Proposition~\ref{prop:positivity} makes the six-dimensional construction finite and reproducible, and the continuation argument shows that truncating the initial data outside the cone does not introduce an earlier singularity. The method is specific to existence of a corotational self-similar solution. It does not establish blowup stability, and it does not address the constant-curvature hyperbolic target or all lower dimensions. In particular, the dimension-six result does not by itself resolve the corresponding questions in dimensions three through five.

\appendix
\section{Exact rational verification}\label{app:code}

For reproducibility, the following Python program verifies every assertion in Table~\ref{tab:cert}. It uses only the standard library, and all polynomial coefficients, endpoints and intermediate arithmetic are exact integers or rational numbers. A polynomial is represented by its list of monomial coefficients in increasing order. The routine \texttt{bernstein} implements \eqref{eq:power}--\eqref{eq:berformula}. Thus the positivity certificate can be checked independently without a numerical ODE solver, a tolerance parameter, or a symbolic algebra package.
\begin{lstlisting}
"""Exact rational positivity certificates for the six-dimensional profile.
Requires Python 3.8+; uses only the standard library.
Polynomials are in increasing monomial order.
"""
from fractions import Fraction as Q
from math import comb


def add(a,b):
    return [(a[i] if i<len(a) else 0)+(b[i] if i<len(b) else 0)
            for i in range(max(len(a),len(b)))]


def scale(a,c):
    return [c*v for v in a]


def mul(a,b):
    out=[Q(0)]*(len(a)+len(b)-1)
    for i,u in enumerate(a):
        for j,v in enumerate(b):
            out[i+j]+=u*v
    return out


def diff(a):
    return [i*a[i] for i in range(1,len(a))]


def bernstein(a,l,r):
    n=len(a)-1
    # coefficients after x=l+(r-l)t
    c=[sum(a[j]*comb(j,k)*l**(j-k)*(r-l)**k
           for j in range(k,n+1)) for k in range(n+1)]
    return [sum(c[k]*Q(comb(i,k),comb(n,k))
                for k in range(i+1)) for i in range(n+1)]

P=[Q(1),Q(7,10),Q(63,4),Q(-1356,5),Q(10836,5),Q(-47502,5),
   Q(248577,10),Q(-199362,5),Q(192483,5),Q(-102807,5),Q(18693,4)]
R=add(mul([0,2,-2],diff(P)),mul([5,-2],P))
P2=mul(P,P)
S=add(mul([1,-1],P2),[4*c/Q(i+1) for i,c in enumerate(P2)])
C=add(mul(add(R,mul([0,2],diff(R))),S),scale(mul(P,mul(R,R)),-1))
assert C[0]==0
D=C[1:]
assert len(D)==32 and D[-1]!=0
for i in range(2):
    b=min(bernstein(P,Q(i,2),Q(i+1,2)))
    assert b>=1
    print('P interval', i, ': exact minimum',b)
for i,L in enumerate([4,13,48,100]):
    b=min(bernstein(D,Q(i,4),Q(i+1,4)))
    assert b>L
    print('D interval', i, ': certified >',L,'; exact minimum',b)
print('D(0) =',D[0])
print('All exact positivity checks passed.')
\end{lstlisting}
The checks give $1$ and $33393/20480$ as the two exact minima for $P$, certify the four strict integer lower bounds for $D$, and return $D(0)=13$. This completes the finite algebraic part of the proof.

\end{document}